# CONFIDENCE BANDS IN NONPARAMETRIC TIME SERIES REGRESSION

By Zhibiao Zhao and Wei Biao Wu[1]

*Pennsylvania State University and University of Chicago*

We consider nonparametric estimation of mean regression and conditional variance (or volatility) functions in nonlinear stochastic regression models. Simultaneous confidence bands are constructed and the coverage probabilities are shown to be asymptotically correct. The imposed dependence structure allows applications in many linear and nonlinear auto-regressive processes. The results are applied to the S&P 500 Index data.

**1. Introduction.** There are two popular approaches in time series analysis: parametric and nonparametric methods. In the literature various parametric models have been proposed, including the classical ARMA, threshold AR (TAR: Tong [32]), exponential AR (EAR: Haggan and Ozaki [21]) and AR with conditional heteroscedasticity (ARCH: Engle [10]) among others. Those models are widely used in practice. An attractive feature of parametric models is that they can provide explanatory insights into the dynamical characteristics of the underlying data-generating mechanism.

However, a parametric model has good performance only when it is indeed the true model or a good approximation of it. Thus, for parametric models, modelling bias may arise and there is a risk of mis-specification that can lead to misunderstanding of the truth and wrong conclusions. One traditional approach is to consider a larger family of parametric models, hoping that bigger models could better approximate the true one. Another appealing way out is to use nonparametric techniques which let the data "speak for themselves" by imposing no specific structures on the underlying regression functions other than smoothness assumptions. See Fan and Yao [17] for an extensive exposition of nonparametric time series analysis.

Received July 2007; revised July 2007.

[1]Supported in part by NSF Grant DMS-04-78704.

*AMS 2000 subject classifications.* Primary 62G08; secondary 62G15.

*Key words and phrases.* Long-range dependence, model validation, moderate deviation, nonlinear time series, nonparametric regression, short-range dependence.

Nonparametric estimates can be used to test parametric assumptions. The problem of nonparametric model validation under dependence is important but difficult. Fan and Yao [17] dealt with this deep problem for time series data by using the idea of generalized likelihood ratio test (Fan, Zhang and Zhang [18]), which is developed for independent data. Fan and Yao [17] pointed out that there have been virtually no theoretical developments on nonparametric model validations under dependence, despite the importance of the latter problem since dependence is an intrinsic characteristic in time series.

Here we consider the model validation problem for the stochastic regression model

$$(1.1) \qquad Y_i = \mu(X_i) + \sigma(X_i)\varepsilon_i, \qquad i = 1, 2, \ldots, n,$$

where $\varepsilon_i$ are independent and identically distributed (i.i.d.) random noises and $(X_i, Y_i)$ are observations. The functions $\mu(\cdot)$ and $\sigma(\cdot)$ are mean regression and conditional variance (or volatility) functions, respectively. As a special case, let $X_i = Y_{i-1}$. Then (1.1) is reduced to the nonlinear autoregressive process $Y_i = \mu(Y_{i-1}) + \sigma(Y_{i-1})\varepsilon_i$ and it includes many parametric time series models. For example, if $\mu(x) = ax$, or $\mu(x) = a\max(x, 0) + b\min(x, 0)$, or $\mu(x) = [a + b\exp(-cx^2)]x$, where $a, b, c$ are real parameters, then it becomes AR, TAR, EAR processes, respectively. For ARCH processes, $\mu = 0$ and $\sigma(x) = (a^2 + b^2x^2)^{1/2}$.

In (1.1), if we let $Y_i = X_{i+1} - X_i$ and $\varepsilon_i$ be i.i.d. standard normals, then (1.1) can be viewed as a discretized version of the stochastic diffusion model

$$(1.2) \qquad dX_t = \mu(X_t)\,dt + \sigma(X_t)\,dW_t,$$

where $\{W_t\}$ is a standard Brownian motion. Many well-known financial models are special cases of (1.2); see Fan [13] and references therein. By model (1.1) we assume that we observe data on an increasing time span $n \to \infty$ with fixed time duration between two consecutive observations. Another approach is sampling from (1.2) on a fixed time span, say $[0, 1]$, with an increasing sampling frequency. These two frameworks have different ranges of applicability. In many applications, we add today's current S&P 500 Index at the end of the existing sample, not the finer observations between two existing ones. In such cases, the former one could be a better choice. See Aït-Sahalia [1] and Bandi and Phillips [3] for more discussions.

We shall address the model validation problem of (1.1) by constructing nonparametric simultaneous confidence bands (SCB) for $\mu$ and $\sigma$. SCB are useful in testing whether $\mu$ and $\sigma$ are of certain parametric forms. For example, in model (1.1), interesting problems include testing whether $\mu$ is linear, quadratic, or of some other pattern and whether $\sigma$ is nonconstant, that is, the existence of conditional heteroscedasticity. The mean regression function $\mu$ can be nonparametrically estimated by kernel, local linear, spline and wavelet methods. To construct asymptotic SCB for $\mu(x)$ over the interval $x \in \mathcal{T} = [T_1, T_2]$ with level $100(1-\alpha)\%$, $\alpha \in (0,1)$, we need to find two functions $l(\cdot) = l_n(\cdot)$ and $u(\cdot) = u_n(\cdot)$ based on the data $(X_i, Y_i)_{i=1}^n$ such that

$$(1.3) \qquad \lim_{n \to \infty} \mathbb{P}\{l(x) \le \mu(x) \le u(x) \text{ for all } x \in \mathcal{T}\} = 1 - \alpha.$$

It is certainly more desirable to have (1.3) in a nonasymptotic sense, namely the probability in (1.3) is exactly $1 - \alpha$. However the latter problem is intractable since it is difficult to establish a finite sample distributional theory for nonparametric regression estimates. With the SCB, we can test whether $\mu$ is of certain parametric form: $H_0 : \mu = \mu_\theta$, where $\theta \in \Theta$ and $\Theta$ is a parameter space. For example, to test whether $\mu(x) = \beta_0 + \beta_1 x$, we can apply the linear regression method and obtain an estimate $(\hat{\beta}_0, \hat{\beta}_1)$ of $(\beta_0, \beta_1)$ from the data $(X_i, Y_i)_{i=1}^n$, and then check whether $l(x) \le \hat{\beta}_0 + \hat{\beta}_1 x \le u(x)$ holds for all $x \in \mathcal{T}$. If so, then we accept at level $\alpha$ the null hypothesis that $\mu$ is linear. Otherwise $H_0$ is rejected.

The construction of SCB $l$ and $u$ satisfying (1.3) has been a difficult problem if dependence is present. Assuming that $(X_i, Y_i)$ are independent random samples from a bivariate population, Johnston [26] obtained an asymptotic distributional theory for $\sup_{0 \le x \le 1} |\hat{\mu}(x) - \mathbb{E}[\hat{\mu}(x)]|$, where $\hat{\mu}(x)$ is the Nadaraya–Watson estimate of the mean regression function $\mu(x) = \mathbb{E}(Y|X=x)$. Johnston applied his limit theorem and constructed asymptotic SCB for $\mu$. Since his result is no longer valid if dependence is present, Johnston's procedure is not applicable in the time series setting. A key tool in Johnston's approach is Bickel and Rosenblatt's [4] asymptotic theory for maximal deviations of kernel density estimators. Bickel and Rosenblatt applied a deep result in probability theory, strong approximation, which asserts that normalized empirical processes of independent random variables can be approximated by Brownian bridges. Such a result generally does not exist under dependence. For other contributions under the independence assumption see Härdle [24], Knafl, Sacks and Ylvisaker [27], Hall and Titterington [23], Härdle and Marron [25], Eubank and Speckman [11], Sun and Loader [31], Xia [38], Cummins, Filloon and Nychka [7] and Dümbgen [9] among others.

In the fixed design case with $X_i = i/n$, Robinson [30] obtained a central limit theorem for nonparametric estimates of trend functions for dependent data. By applying Komlós et al.'s [28] strong invariance principle for partial sums, Eubank and Speckman [11] constructed SCB for $\mu$ with asymptotically correct coverage probabilities. Their method was extended to the time series setting by Wu and Zhao [37]. However, Wu and Zhao's result is not applicable here since it heavily relies on the fixed design assumption.

In this paper, we shall consider a variant of (1.3) and construct SCB over a subset $\mathcal{T}_n$ of $\mathcal{T}$ with $\mathcal{T}_n$ becoming denser as $n \to \infty$. A similar framework is also adopted in Bühlmann [5] and Knafl, Sacks and Ylvisaker [27] in other contexts. The latter paper concerns the construction of conservative SCB of regression functions for independent data with Gaussian errors over grid points. It is shown that our SCB has asymptotically correct coverage probabilities under a general dependence structure on $(X_i, Y_i)$ which allows applications in many linear and nonlinear processes. Our method can be applied to statistical inference problems in time series including goodness-of-fit, hypothesis testing and others. In the development of our asymptotic theory, we apply a deep martingale moderate deviation principle by Grama and Haeusler [19].

We now introduce some notation. For a random variable $Z$ write $Z \in \mathcal{L}^p, p > 0$, if $\|Z\|_p := [\mathbb{E}(|Z|^p)]^{1/p} < \infty$, and $\|Z\| = \|Z\|_2$. For $a, b \in \mathbb{R}$ let $a \wedge b = \min(a, b), a \vee b = \max(a, b)$ and $\lceil a \rceil = \inf\{k \in \mathbb{Z} : k \geq a\}$. Let $\{a_n\}$ and $\{b_n\}$ be two real sequences. We write $a_n \asymp b_n$ if $|a_n/b_n|$ is bounded away from 0 and $\infty$ for all large $n$. For $\mathcal{S} \subset \mathbb{R}$ denote by $\mathcal{C}^p(\mathcal{S}) = \{g(\cdot) : \sup_{x \in \mathcal{S}} |g^{(k)}(x)| < \infty, k = 0, \ldots, p\}$ the set of functions having bounded derivatives on $\mathcal{S}$ up to order $p \geq 1$, and by $\mathcal{C}^0(\mathcal{S})$ the set of continuous functions on $\mathcal{S}$. Let $\mathcal{S}^\epsilon = \bigcup_{y \in \mathcal{S}}\{x : |x - y| \leq \epsilon\}$ be the $\epsilon$-neighborhood of $\mathcal{S}, \epsilon > 0$.

The rest of the paper is structured as follows: We introduce our dependence structure on $(X_i, Y_i, \varepsilon_i)$ in Section 2. Section 3 presents the main results: SCB for $\mu(\cdot)$ and $\sigma^2(\cdot)$ with asymptotically correct coverage probabilities are constructed in Sections 3.1 and 3.2, respectively. In Section 4, applications are made to two important cases of (1.1): nonlinear time series and linear processes, where we consider both short- and long-range dependent processes. In Section 5, we discuss some implementation issues, and then perform a simulation study. Section 6 contains an application to the S&P 500 Index data. We defer the proofs to Section 7.

**2. Dependence structure.** In (1.1), assume that $\varepsilon_i, i \in \mathbb{Z}$, are i.i.d. and that $X_i$ is a stationary process

$$(2.1) \qquad X_i = G(\mathcal{F}_i) \qquad \text{where } \mathcal{F}_i = (\ldots, \eta_{i-1}, \eta_i).$$

Here $\eta_i, i \in \mathbb{Z}$, are i.i.d. and $G$ is a measurable function such that $X_i$ is well defined. The framework (2.1) is very general (Priestley [29], Tong [32], Wu [35]). Assume that $\varepsilon_i$ is independent of $\mathcal{F}_i$ and $\eta_i$ is independent of $\varepsilon_j, j \leq i - 2$.

Define the projection operator $\mathcal{P}_k, k \in \mathbb{Z}$, by $\mathcal{P}_k Z = \mathbb{E}(Z|\mathcal{F}_k) - \mathbb{E}(Z|\mathcal{F}_{k-1})$, $Z \in \mathcal{L}^1$. Let $F_X$ and $F_\varepsilon$ be the distribution functions of $X_0$ and $\varepsilon_0$, respectively; let $f_X = F_X'$ and $f_\varepsilon = F_\varepsilon'$ be their densities. Let $F_X(x|\mathcal{F}_i) = \mathbb{P}(X_{i+1} \leq x|\mathcal{F}_i), i \in \mathbb{Z}$, be the conditional distribution function of $X_{i+1}$ given $\mathcal{F}_i$ and $f_X(x|\mathcal{F}_i) = \partial F_X(x|\mathcal{F}_i)/\partial x$ the conditional density. Let

$$(2.2) \qquad \theta_i = \sup_{x \in \mathbb{R}} \|\mathcal{P}_0 f_X(x|\mathcal{F}_i)\| + \sup_{x \in \mathbb{R}} \|\mathcal{P}_0 f_X'(x|\mathcal{F}_i)\|,$$

where $f'_X(x|\mathcal{F}_i) = \partial f_X(x|\mathcal{F}_i)/\partial x$ if it exists. For $n \in \mathbb{N}$ define

$$(2.3) \qquad \Xi_n = n\Theta_{2n}^2 + \sum_{k=n}^{\infty}(\Theta_{n+k} - \Theta_k)^2 \qquad \text{where } \Theta_n = \sum_{i=1}^{n}\theta_i.$$

Roughly speaking, $\theta_i$ measures the contribution of $\varepsilon_0$ in predicting $X_{i+1}$ (Wu [35]). If $\Theta_\infty < \infty$, then the cumulative contribution of $\varepsilon_0$ in predicting future values is finite, thus implying short-range dependence (SRD). In this case $\Xi_n = O(n)$. Our setting also allows long-range dependence (LRD). For example, let $\theta_i = i^{-\beta}\ell(i)$, where $\beta > 1/2$ and $\ell(\cdot)$ is a slowly varying function, namely $\lim_{x \to \infty} \ell(\lambda x)/\ell(x) = 1$ for all $\lambda > 0$. Note that $\bar{\bar{\ell}}(n) = \sum_{i=1}^{n}|\ell(i)|/i$ is also a slowly varying function. By Karamata's theorem,

$$(2.4) \qquad \Xi_n = O(n), \qquad O[n^{3-2\beta}\ell^2(n)] \quad \text{or} \quad O\{n[\bar{\bar{\ell}}(n)]^2\},$$

under $\beta > 1$ (SRD case), $\beta < 1$ (LRD case) or $\beta = 1$, respectively (see Wu [34]). In the LRD case $\Xi_n$ grows faster than $n$. In Section 4 we shall give bounds on $\Xi_n$ for SRD and LRD linear processes and some nonlinear time series.

**3. Main results.** Let $\mathcal{T} = [T_1, T_2]$ be a bounded interval. We assume hereafter without loss of generality (WLOG) that $\mathbb{E}(\varepsilon_0) = 0$ and $\mathbb{E}(\varepsilon_0^2) = 1$ since otherwise model (1.1) can be re-parameterized by letting $\bar{\mu}(x) = \mu(x) + \sigma(x)\mathbb{E}(\varepsilon_0)$, $\bar{\sigma}(x) = c\sigma(x)$ and $\bar{\varepsilon}_i = [\varepsilon_i - \mathbb{E}(\varepsilon_i)]/c$, where $c^2 = \mathbb{E}(\varepsilon_0^2) - [\mathbb{E}(\varepsilon_0)]^2$.

3.1. *Simultaneous confidence band for $\mu$.* There exists a vast literature on nonparametric estimation of the regression function $\mu$. Here we use the Nadaraya–Watson estimator

$$(3.1) \qquad \hat{\mu}_{b_n}(x) = \frac{1}{nb_n\hat{f}_X(x)}\sum_{i=1}^{n}K_{b_n}(x - X_i)Y_i$$

$$\text{where } \hat{f}_X(x) = \frac{1}{nb_n}\sum_{i=1}^{n}K_{b_n}(x - X_i).$$

Here and hereafter $K_{b_n}(u) = K(u/b_n)$, $K$ is a kernel function with $\int_{\mathbb{R}} K(u)\,du = 1$ and the bandwidth $b_n \to 0$ satisfies $nb_n \to \infty$. One can also use the local linear estimator in Fan and Gijbels [15] and derive similar results. In this paper we have decided to use the local constant estimator (3.1) instead of the local linear estimator since the latter involves more tedious theoretical derivations. In Definition 1 below, some regularity conditions on $K$ are imposed. Theorem 1 asserts a central limit theorem (CLT) for $\hat{\mu}_{b_n}(x)$, which can be used to construct point-wise confidence intervals for $\mu(x)$.

DEFINITION 1. *Let $\mathcal{K}$ be the set of kernels which are bounded, symmetric, and have bounded derivative and bounded support. Let $\psi_K = \int_{\mathbb{R}} u^2 K(u)\, du/2$ and $\varphi_K = \int_{\mathbb{R}} K^2(u)\, du$.*

THEOREM 1. *Let $x \in \mathbb{R}$ be fixed and $K \in \mathcal{K}$. Assume that $f_X(x) > 0, \sigma(x) > 0$ and $f_X, \mu \in \mathcal{C}^4(\{x\}^\epsilon)$ for some $\epsilon > 0$. Further assume that*

$$(3.2) \qquad nb_n^9 + \frac{1}{nb_n} + \Xi_n\left(\frac{b_n^3}{n} + \frac{1}{n^2}\right) \to 0.$$

*Let $\rho_\mu(x) = \mu''(x) + 2\mu'(x) f_X'(x)/f_X(x)$. Then as $n \to \infty$,*

$$\frac{\sqrt{nb_n}}{\sigma(x)\sqrt{\varphi_K}}\sqrt{\hat{f}_X(x)}[\hat{\mu}_{b_n}(x) - \mu(x) - b_n^2\psi_K\rho_\mu(x)] \Rightarrow N(0,1).$$

THEOREM 2. *Let $\varepsilon_0 \in \mathcal{L}^3$, $\mathcal{T} = [T_1, T_2]$ and $K \in \mathcal{K}$. Assume that $\inf_{x \in \mathcal{T}} f_X(x) > 0, \inf_{x \in \mathcal{T}} \sigma(x) > 0$, and $f_X, \mu \in \mathcal{C}^4(\mathcal{T}^\epsilon), \sigma \in \mathcal{C}^2(\mathcal{T}^\epsilon)$ for some $\epsilon > 0$. Further assume*

$$(3.3) \qquad nb_n^9 \log n + \frac{(\log n)^3}{nb_n^3} + \Xi_n\left[\frac{b_n^3 \log n}{n} + \frac{(\log n)^2}{n^2 b_n^{4/3}}\right] \to 0.$$

*Let $\rho_\mu(x)$ be as in Theorem 1. For $n \geq 2$ define*

$$(3.4) \quad B_n(z) = \sqrt{2\log n} - \frac{1}{\sqrt{2\log n}}\left[\frac{1}{2}\log\log n + \log(2\sqrt{\pi})\right] + \frac{z}{\sqrt{2\log n}}.$$

*Let the kernel $K$ have support $[-k_0, k_0]$; let $\mathcal{T}_n = \{x_j = T_1 + 2k_0 b_n j, j = 0, 1, \ldots, m_n - 1$ and $m_n = \lceil (T_2 - T_1)/(2k_0 b_n) \rceil$. Then for every $z \in \mathbb{R}$,*

$$\lim_{n\to\infty} \mathbb{P}\left\{\frac{\sqrt{nb_n}}{\sqrt{\varphi_K}} \sup_{x \in \mathcal{T}_n} \frac{[\hat{f}_X(x)]^{1/2}}{\sigma(x)}|\hat{\mu}_{b_n}(x) - \mu(x) - b_n^2\psi_K\rho_\mu(x)| \leq B_{m_n}(z)\right\}$$
$$= e^{-2e^{-z}}.$$

Observe that $\mathcal{T}_n$ becomes denser in $\mathcal{T}$ as $b_n \to 0$. Since $b_n \to 0$, if the regression function $\mu$ is sufficiently smooth, then $\{\mu(x) : x \in \mathcal{T}\}$ can be well approximated by $\{\mu(x) : x \in \mathcal{T}_n\}$ for large $n$. Theorem 2 is useful to construct SCB in an approximate version of (1.3):

$$(3.5) \qquad \lim_{n\to\infty} \mathbb{P}\{l(x) \leq \mu(x) \leq u(x) \text{ for all } x \in \mathcal{T}_n\} = 1 - \alpha.$$

Specifically, let $\hat{\sigma}_n(x)$ [resp. $\hat{\rho}_\mu(x)$] be an estimate of $\sigma(x)$ [resp. $\rho_\mu(x)$] such that $\sup_{x \in \mathcal{T}} |\hat{\sigma}_n(x) - \sigma(x)| = o_p[(\log n)^{-1}]$ and $\sup_{x \in \mathcal{T}} |\hat{\rho}_\mu(x) - \rho_\mu(x)| = o_p[(nb_n^5 \log n)^{-1/2}]$. By Slutsky's theorem, Theorem 2 still holds if $\sigma$ and $\rho_\mu$

therein is replaced by $\hat{\sigma}_n$ and $\hat{\rho}_\mu$, respectively. Hence, in the sense of (3.5), an asymptotic $100(1-\alpha)\%$ SCB for $\mu$ can be constructed as

$$
\hat{\mu}_{b_n}(x) - b_n^2 \psi_K \hat{\rho}_\mu(x) \pm \frac{\sqrt{\varphi_K} \hat{\sigma}_n(x)}{\sqrt{n b_n \hat{f}_X(x)}} B_{m_n}(z_\alpha) \quad \text{and}
$$

(3.6)

$$
z_\alpha = -\log\log[(1-\alpha)^{-1/2}].
$$

In (3.6), $\rho_\mu(x)$ cannot be easily estimated since it involves unknown functions $\mu''$, $\mu'$ and $f_X'$. Following Wu and Zhao [37], we adopt the simple jackknife-type bias correction procedure which avoids estimating $\mu''$, $\mu'$ and $f_X'$:

(3.7)
$$
\hat{\mu}_{b_n}^*(x) = 2\hat{\mu}_{b_n}(x) - \hat{\mu}_{\sqrt{2}b_n}(x).
$$

Then the bias term $O(b_n^2)$ in $\hat{\mu}_{b_n}$ reduces to $O(b_n^4)$ in $\hat{\mu}_{b_n}^*$. Using (3.7) is equivalent to using the 4th-order kernel $K^*(u) = 2K(u) - K(u/\sqrt{2})/\sqrt{2}$ in (3.1). Clearly, $K^* \in \mathcal{K}$ has support $[-\sqrt{2}k_0, \sqrt{2}k_0]$ and $\psi_{K^*} = 0$. Let $m_n^* = \lceil (T_2 - T_1)/(2\sqrt{2}k_0 b_n) \rceil$ and $\mathcal{T}_n^* = \{x_j^* = T_1 + 2\sqrt{2}k_0 b_n j, j = 0, 1, \ldots, m_n^* - 1\}$. Then Theorem 2 still holds with $\hat{\mu}_{b_n}$ (resp. $K, m_n, \mathcal{T}_n$) replaced by $\hat{\mu}_{b_n}^*$ (resp. $K^*, m_n^*, \mathcal{T}_n^*$).

In Theorem 2, (3.3) imposes conditions on the bandwidth $b_n$ and the strength of the dependence. The first part $n b_n^9 \log n \to 0$ aims to control the bias with $b_n$ being not too large, while the second one $(\log n)^3/(n b_n^3) \to 0$ suggests that $b_n$ should not be too small, thus ensuring the validity of the moderate deviation principle (see the proof of Theorem 5). The third part suggests that the dependence should not be too strong. For SRD processes, we have $\Xi_n = O(n)$ and the third term in (3.3) is automatically $o(1)$ if the first two are $o(1)$. If $b_n \asymp n^{-\beta}$ with $\beta \in (1/9, 1/3)$, then the first two terms in (3.3) are $o(1)$. In particular, (3.3) allows $\beta = 1/5$, which corresponds to the mean square error (MSE)-optimal bandwidth. Interestingly, (3.3) also allows long-range dependent processes; see Section 4.2.

3.2. *Simultaneous confidence band for $\sigma^2$.* Let $\hat{\mu}_{b_n}^*$ be as in (3.7). Since $\mathbb{E}(\varepsilon_i^2) = 1$ and $\mathbb{E}\{[Y_i - \mu(X_i)]^2 | X_i = x\} = \sigma^2(x)$, a natural residual-based estimator of $\sigma^2(x)$ is

$$
\hat{\sigma}_{h_n}^2(x) = \frac{1}{n h_n \tilde{f}_X(x)} \sum_{i=1}^n [Y_i - \hat{\mu}_{b_n}^*(X_i)]^2 K_{h_n}(x - X_i),
$$

(3.8)

$$
\text{where } \tilde{f}_X(x) = \frac{1}{n h_n} \sum_{i=1}^n K_{h_n}(x - X_i).
$$

Here $h_n$ is another bandwidth and it can be different from the bandwidth $b_n$ in estimating $\mu$. Similar methods are applied in Fan and Yao [16] and Hall and Carroll [22].

REMARK 1. A referee pointed out that, as shown in Fan and Yao [16], it is not necessary to correct the bias in estimating the mean function $\mu$. In fact, if one uses the original estimate $\hat{\mu}_{b_n}$, then the bias term $O(b_n^2)$ on the right-hand side of (7.23) is squared in (7.24) and becomes $O(b_n^4)$, which does not affect the subsequent proofs.

PROPOSITION 1. Let $K \in \mathcal{K}, \varepsilon_0 \in \mathcal{L}^6$ and $h_n \asymp b_n$. Assume that $\inf_{x \in \mathcal{T}} f_X(x) > 0$, $\inf_{x \in \mathcal{T}} \sigma(x) > 0$ and $f_X, \mu \in \mathcal{C}^4(\mathcal{T}^\epsilon)$ for some $\epsilon > 0$. Further assume that

$$(3.9) \qquad h_n^{3/2} \log n + \frac{1}{n^2 h_n^5} + \frac{\Xi_n}{n^2} \to 0.$$

Then

$$\sup_{x \in \mathcal{T}} |\hat{\sigma}_{h_n}^2(x) - \sigma^2(x)| = O_p \left\{ h_n^2 + \frac{1}{nh_n^{5/2}} + \left[ \frac{\log n}{nh_n} \right]^{1/2} + \left[ \frac{\log n}{n^3 h_n^7} \right]^{1/4} + \frac{\Xi_n^{1/2} h_n}{n} \right\}.$$

Proposition 1 provides a uniform error bound for the estimate $\hat{\sigma}_{h_n}^2(\cdot)$. From the proof of Proposition 1 and Theorem 3, one can as in Theorem 1 establish a CLT for $\hat{\sigma}_{h_n}^2(x)$ for each fixed $x$ and the optimal bandwidth $h_n \asymp n^{-1/5}$. We omit the details. In Proposition 1, if one uses the optimal bandwidth $h_n \asymp n^{-1/5}$, then $\sup_{x \in \mathcal{T}} |\hat{\sigma}_{h_n}^2(x) - \sigma^2(x)| = O_p[n^{-2/5}(\log n)^{1/2} + \Xi_n^{1/2} n^{-6/5}]$. The first part $O_p[n^{-2/5}(\log n)^{1/2}]$ in the error bound is optimal in nonparametric curve estimation for independent data. The second part accounts for dependence, and it can be absorbed into the first one if $\Xi_n = O(n^{8/5})$. In particular, for SRD processes, $\Xi_n = O(n)$.

Theorem 3 below presents a maximal deviation result for $\hat{\sigma}_{h_n}^2$ and it can be used to construct SCB for $\sigma^2$.

THEOREM 3. Let the conditions in Proposition 1 be fulfilled. Further assume that $\sigma \in \mathcal{C}^4(\mathcal{T}^\epsilon)$ for some $\epsilon > 0$ and

$$(3.10) \qquad nh_n^9 \log n + \frac{\log n}{nh_n^4} + \Xi_n \left[ \frac{h_n^3 \log n}{n} + \frac{(\log n)^2}{n^2 h_n^{4/3}} \right] \to 0.$$

Let the kernel $K$ have support $[-k_0, k_0]$; let $\tilde{\mathcal{T}}_n = \{\tilde{x}_j = T_1 + 2k_0 h_n j, j = 0, 1, \dots, \tilde{m}_n - 1\}$ and $\tilde{m}_n = \lceil (T_2 - T_1)/(2k_0 h_n) \rceil$. Let $B_n(z)$ be as in (3.4). Then for every $z \in \mathbb{R}$,

$$\lim_{n \to \infty} \mathbb{P} \left\{ \frac{\sqrt{nh_n}}{\sqrt{\varphi_K \nu_\varepsilon}} \sup_{x \in \tilde{\mathcal{T}}_n} \frac{[\tilde{f}_X(x)]^{1/2}}{\hat{\sigma}_{h_n}^2(x)} |\hat{\sigma}_{h_n}^2(x) - \sigma^2(x) - h_n^2 \psi_K \rho_\sigma(x)| \le B_{\tilde{m}_n}(z) \right\}$$
$$= e^{-2e^{-z}},$$

*where $\nu_\varepsilon = \mathbb{E}(\varepsilon_0^4) - 1 > 0$ and*

$$\rho_\sigma(x) = 2\sigma'(x)^2 + 2\sigma(x)\sigma''(x) + 4\sigma(x)\sigma'(x)f_X'(x)/f_X(x).$$

If $\mu$ were known and we use $\hat{\sigma}_{h_n}^2$ in (3.8) to estimate $\sigma^2$ with $\hat{\mu}_{b_n}$ therein replaced by the true function $\mu$, then Theorem 3 is still applicable. So Theorem 3 implies the oracle property that the construction of SCB for $\sigma^2$ does not heavily rely on the estimation of $\mu$. Under strong mixing conditions, Fan and Yao [16] obtained a similar oracle property for their local linear estimator of $\sigma^2$. For independent data, Hall and Carroll [22] considered the effect of estimating the mean on variance function estimation.

As in (3.7), we propose the bias-corrected estimate $\hat{\sigma}_{h_n}^{2*}(x) = 2\hat{\sigma}_{h_n}^2(x) - \hat{\sigma}_{\sqrt{2}h_n}^2$. Similarly as in Section 3.1, we can define $\tilde{m}_n^*$ and $\tilde{\mathcal{T}}_n^*$ accordingly and Theorem 3 still holds with $\hat{\sigma}_{h_n}^2$ (resp. $K, \tilde{m}_n, \tilde{\mathcal{T}}_n$) replaced by $\hat{\sigma}_{h_n}^{2*}$ (resp. $K^*, \tilde{m}_n^*, \tilde{\mathcal{T}}_n^*$).

3.3. *Estimation of $\nu_\varepsilon$ in Theorem 3.* To apply Theorem 3, one needs to estimate $\nu_\varepsilon = \mathbb{E}(\varepsilon_0^4) - 1$. Here we estimate $\nu_\varepsilon$ by

$$(3.11) \qquad \hat{\nu}_\varepsilon = \frac{\sum_{i=1}^n \hat{\varepsilon}_i^4 \mathbf{1}_{X_i \in \mathcal{T}}}{\sum_{i=1}^n \mathbf{1}_{X_i \in \mathcal{T}}} - 1$$

$$\text{where } \hat{\varepsilon}_i = \frac{Y_i - \hat{\mu}_{b_n}^*(X_i)}{\hat{\sigma}_{h_n}^*(X_i)}, \ i = 1, 2, \ldots, n.$$

Here $\hat{\varepsilon}_i$ are estimated residuals for model (1.1). The naive estimate $n^{-1}\sum_{i=1}^n \hat{\varepsilon}_i^4 - 1$ does not have a good practical performance since $\hat{\sigma}_{h_n}^*(x)$ behaves poorly if $|x|$ is too large. Truncation by $\mathcal{T}$ improves the performance.

PROPOSITION 2. *Assume that the conditions in Proposition 1 are satisfied. Then*

$$(3.12)\hat{\nu}_\varepsilon - \nu_\varepsilon = O_\mathrm{p}\left\{ n^{-1/3} + h_n^4 + \frac{1}{nh_n^{5/2}} + \left[\frac{\log n}{nh_n}\right]^{1/2} + \left[\frac{\log n}{n^3 h_n^7}\right]^{1/4} + \frac{\Xi_n^{1/2}}{n} \right\}.$$

By Proposition 2, when one chooses the MSE-optimal bandwidths $b_n \asymp h_n \asymp n^{-1/5}$ and assume $\Xi_n = O(n^{4/3})$, then $\hat{\nu}_\varepsilon - \nu_\varepsilon = O_\mathrm{p}(n^{-1/3})$. By Slutsky's theorem, the convergence in Theorem 3 still holds when $\nu_\varepsilon$ therein is replaced by $\hat{\nu}_\varepsilon$. It is expected that, under stronger regularity conditions, one can obtain the better bound $\hat{\nu}_\varepsilon - \nu_\varepsilon = O_\mathrm{p}(n^{-p})$ with $p > 1/3$. Here we do not pursue this line of direction since the involved calculation will be tedious and since the bound $\hat{\nu}_\varepsilon - \nu_\varepsilon = O_\mathrm{p}[(\log n)^{-1}]$ suffices for our purpose.

**4. Examples.** To apply Theorems [2] and [3], we need to deal with $\Xi_n$ defined in (2.3). Let $\eta'_0, \eta_i, i \in \mathbb{Z}$, be i.i.d.; let $\mathcal{F}'_i = (\mathcal{F}_{-1}, \eta'_0, \eta_1, \ldots, \eta_i)$. By Theorem 1 in Wu [35], we have

$$(4.1) \quad \theta_i \le \varpi_i := \sup_{x \in \mathbb{R}} \|f_X(x|\mathcal{F}_i) - f_X(x|\mathcal{F}'_i)\| + \sup_{x \in \mathbb{R}} \|f'_X(x|\mathcal{F}_i) - f'_X(x|\mathcal{F}'_i)\|.$$

For many processes there exist simple and easy-to-use bounds for $\varpi_i$. Here we shall consider linear processes and some popular nonlinear time series models.

4.1. *Short-range dependent linear processes.* Let $\eta_i, i \in \mathbb{Z}$, be i.i.d. with $\eta_0 \in \mathcal{L}^q$, $q > 0$, and $\mathbb{E}(\eta_0) = 0$ if $q \ge 1$. For real sequence $(a_i)_{i \ge 0}$ satisfying $\sum_{i=0}^\infty |a_i|^{q \wedge 2} < \infty$, the process

$$(4.2) \qquad\qquad X_i = \sum_{j=0}^\infty a_j \eta_{i-j},$$

is well defined and stationary. Special cases of (4.2) include ARMA models. Assume WLOG that $a_0 = 1$. Let $\bar{X}_i = X_i - \eta_i$ and $\bar{X}'_i = \bar{X}_i + a_i(\eta'_0 - \eta_0)$. Then $f_X(x|\mathcal{F}_{i-1}) = f_\eta(x - \bar{X}_i)$ and $f_X(x|\mathcal{F}'_{i-1}) = f_\eta(x - \bar{X}'_i)$, where $f_\eta$ is the density of $\eta_0$. Assume that $f_\eta \in \mathcal{C}^2(\mathbb{R})$. Simple calculations show that $\theta_i = O(|a_i|^{q'})$, where $q' = (q \wedge 2)/2$; see Proposition 2 in Zhao and Wu [39]. Therefore we have $\Xi_n = O(n)$ if $\sum_{i=1}^\infty |a_i|^{q'} < \infty$. If $q \ge 2$, then the latter condition becomes $\sum_{i=0}^\infty |a_i| < \infty$. For causal ARMA processes, $a_i \to 0$ geometrically quickly. Note that our setting allows heavy-tailed $\eta_i$.

4.2. *Long-range dependent linear processes.* Consider the linear process (4.2) with $a_i = i^{-\alpha} \ell(i)$, where $\alpha > 1/(2q')$, $q' = (q \wedge 2)/2$, and $\ell(\cdot)$ is a slowly varying function. The case of $\alpha q' > 1$ is covered by Section 4.1. Assume $\alpha q' \in (1/2, 1]$. If $q \ge 2$ and $\alpha \in (1/2, 1)$, by Karamata's theorem, the covariances $\mathbb{E}(X_0 X_n)$ are of order $n^{1-2\alpha} \ell^2(n)$ and not summable, hence $(X_i)$ is long-range dependent. As in Section 4.1, $\theta_i = O[i^{-\alpha q'} \ell^{q'}(i)]$. By (2.4), $\Xi_n = O[n^{3-2\alpha q'} \ell^{2q'}(n)]$ if $\alpha q' \in (1/2, 1)$ and $\Xi_n = O\{n[\sum_{i=1}^n |\ell^{q'}(i)|/i]^2\}$ if $\alpha q' = 1$.

If $\alpha q' \in (17/26, 1]$, then (3.3) and (3.10) hold if $b_n \asymp h_n \asymp n^{-\beta}$ and $\beta$ satisfies

$$(4.3) \qquad \max\left\{\frac{1}{9}, \frac{2(1-\alpha q')}{3}\right\} < \beta < \min\left\{\frac{1}{3}, \frac{3(2\alpha q' - 1)}{4}\right\}.$$

Since $\alpha q' \in (17/26, 1]$, such $\beta$ exists. So under (4.3), Theorems [2] and [3] are applicable.

EXAMPLE 1. Let $\alpha(z) = 1 - \sum_{i=1}^k \alpha_i z^i$ and $\beta(z) = 1 + \sum_{i=1}^p \beta_i z^i$ be two polynomials, where $\alpha_1, \ldots, \alpha_k, \beta_1, \ldots, \beta_p \in \mathbb{R}$, are real coefficients. Denote by $B$ the backward shift operator: $B^j X_n = X_{n-j}, j \ge 0$. Consider the FARIMA$(k, d, p)$ process $X_n$ defined by $\alpha(B)(1 - B)^d X_n = \beta(B)\varepsilon_n$, $d \in (-1/2, 1/2)$. Let $\Gamma(x) = \int_0^\infty t^{x-1} e^{-t} \, dt$ be the gamma function. In the simple case of $p = k = 0$, we have $X_n = \sum_{i=0}^\infty a_i \varepsilon_{n-i}$, where

$$(4.4) \qquad a_n = \frac{\Gamma(n + d)}{\Gamma(n + 1)\Gamma(d)} \asymp n^{d-1}.$$

If $d \in (0, 1/2)$, then $X_n$ is long-range dependent. More generally, it can be shown that (4.4) holds for general FARIMA$(k, d, p)$ processes if $\alpha(z) \neq 0$ for all complex $|z| \leq 1$.

4.3. *Nonlinear AR models.* Consider a special case of (1.1) with $X_i = Y_{i-1}$ and $\eta_i = \varepsilon_{i-1}$. Then it becomes the nonlinear AR model

$$(4.5) \qquad X_{i+1} = \mu(X_i) + \sigma(X_i)\eta_{i+1}.$$

Special cases of (4.5) include linear AR, ARCH, TAR and EAR processes. Denote by $f_\eta$ the density of $\eta_0$. Assume $\eta_0 \in \mathcal{L}^q$ and $\sup_{x \in \mathbb{R}}(1 + |x|)[|f_\eta'(x)| + |f_\eta''(x)|] < \infty$. As in Zhao and Wu [39], we have $\theta_i = O(r^i)$ with $r \in (0, 1)$, and hence $\Xi_n = O(n)$, provided that

$$(4.6) \qquad \begin{aligned} &\inf_{x \in \mathbb{R}} \sigma(x) > 0, \qquad \sup_{x \in \mathbb{R}}[|\mu'(x)| + |\sigma'(x)|] < \infty \quad \text{and} \\ &\sup_{x \in \mathbb{R}} \|\mu'(x) + \sigma'(x)\eta_0\|_q < 1. \end{aligned}$$

EXAMPLE 2. Consider the ARCH model $X_n = \eta_n \sqrt{a^2 + b^2 X_{n-1}^2}$, where $\eta_i, i \in \mathbb{Z}$, are i.i.d. and $a, b$ are real parameters. If $\eta_0 \in \mathcal{L}^q$ and $|b| \|\eta_0\|_q < 1$, then (4.6) holds.

5. **A simulation study.** In this section we shall present a simulation study for the performance of our SCB constructed in Section 3. Let $\varepsilon_i, i \in \mathbb{Z}$, be i.i.d. standard normal random variables. We shall consider the following two models:

Model 1: AR(1)    $Y_i = \mu(Y_{i-1}) + s\varepsilon_i$,      $i = 1, 2, \ldots, n$.

Model 2: ARCH(1)    $Y_i = \sigma(Y_{i-1})\varepsilon_i$,      $i = 1, 2, \ldots, n$.

Here $\mu$ and $\sigma$ are functions of interest and $s > 0$ is the scale parameter. Model 1 is a nonlinear AR model and Model 2 is an ARCH model.

The jackknife bias-correction scheme reduces bias and allows one to choose a relatively larger bandwidth. On the other hand, a larger bandwidth results in relatively fewer grid points in $\mathcal{T}_n$, and consequently a less accurate approximation of $\{\mu(x) : x \in \mathcal{T}\}$ by $\{\mu(x) : x \in \mathcal{T}_n\}$. In our simulation, we tried different bandwidths and different sets $\mathcal{T}_{(k)}, k = 20, 30$ and $50$, of grid points to access the performance of our SCB. Here $\mathcal{T}_{(k)}$ denotes the set containing $k$ grid points evenly spaced over $\mathcal{T}$, regardless of the bandwidth. The three choices of $\mathcal{T}_{(k)}$ have similar performance, and the results reported below are for $\mathcal{T}_{(20)}$.

To construct SCB based on Theorem 2, one can use the cutoff value $q_\alpha = B_{m_n}(z_\alpha)$, where $z_\alpha = -\log[(1-\alpha)^{-1/2}]$ is the $(1-\alpha)$th quantile of the limiting extreme value distribution. Since the convergence to extreme value distribution is quite slow, we shall propose a finite sample approximation scheme to compute the cutoff value. Let $Z_i, 1 \le i \le m_n$, be i.i.d. standard normal random variables and $B_{m_n}(z)$ be as in Theorem 2. Then

$$(5.1) \qquad \lim_{m_n \to \infty} \mathbb{P}\left\{ \sup_{1 \le i \le m_n} |Z_i| \le B_{m_n}(z) \right\} = e^{-2e^{-z}},$$

and the events in (3.5) and (5.1) have the same limiting distribution. So we propose the finite sample cutoff value $q_\alpha^*$ defined by

$$(5.2) \quad \mathbb{P}\left\{ \sup_{1 \le i \le m_n} |Z_i| \le q_\alpha^* \right\} = 1 - \alpha, \quad \text{or} \quad \mathbb{P}\{|Z_1| \le q_\alpha^*\} = (1-\alpha)^{1/m_n}.$$

The difference between $q_\alpha$ and $q_\alpha^*$ is that $q_\alpha^*$ is the $(1-\alpha)$th quantile of $\sup_{1 \le i \le m_n} |Z_i|$ for fixed $m_n$ while $q_\alpha$ is based on the limiting distribution. Thus, we expect the SCB based on $q_\alpha^*$ would outperform the one based on $q_\alpha$.

To assess the performance of our SCB, we generate $10^4$ realizations from Model 1 with $n = 2500, \mu(x) = 0.9\sin(x)$ and $s = 0.4$. Then (4.6) holds. Under this setting, simulations show that about 92–95% of the $Y$'s lie within the interval $[-1.1, 1.1]$. Thus we take $\mathcal{T} = [-1.1, 1.1]$. We construct a SCB for each realization, and simulated coverage probabilities are the proportion of these $10^4$ SCBs that cover $\mu$. When applying the simulation procedure, we adopt the following technique for fitting $\hat{\mu}_{b_n}^*$ at $Y_i, 0 \le i \le n$: we fit 300 grid points evenly spaced on the range of $Y_i$'s and use the fitted value of the nearest grid point to $Y_i$ as $\hat{\mu}_{b_n}^*(Y_i)$. Doing this allows one to gain better smoothness since the original series $(Y_i)_{i=0}^n$ may be irregularly spaced. The result is reported in Table 1. The first row corresponds to different choices of bandwidth $b_n$, the second and third rows are the simulated coverage probabilities of the constructed SCBs by using the cutoff values $q_\alpha$ and $q_\alpha^*$, respectively. For $\alpha = 0.05$ and $m_n = 20$, $q_{0.05} = 3.203$ and $q_{0.05}^* = 3.016$. The coverage probabilities using $q_\alpha^*$ are relatively insensitive to the choice of bandwidths and very close to the nominal level 95% while the coverage probabilities using $q_\alpha$ are systematically bigger.

In Model 2, we take $n = 2500$, $\sigma(x) = (0.4 + 0.2x^2)^{1/2}$ and $\mathcal{T} = [-1, 1]$. Table 2 shows similar phenomena.

TABLE 1
*Simulated coverage probabilities for Model 1*

| $b_n$ | 0.10 | 0.12 | 0.14 | 0.15 | 0.16 | 0.18 | 0.20 |
|---|---|---|---|---|---|---|---|
| use $q_\alpha$ | 0.9718 | 0.9716 | 0.9705 | 0.9695 | 0.9671 | 0.9646 | 0.9581 |
| use $q_\alpha^*$ | 0.9471 | 0.9498 | 0.9482 | 0.9479 | 0.9463 | 0.9430 | 0.9312 |

**6. Application to the S&P 500 Index data.** The dataset contains 14374 records, $S_0$, $S_1$, ..., $S_{14373}$, of S&P 500 Index daily data during the period 3 January 1950 to 20 February 2007. Let $Y_i = \log(S_{i+1}) - \log(S_i)$, $i = 0, 1, \ldots, 14372$, be the log returns (difference of the logarithm of prices). Ding, Granger and Engle [8] considered daily log return series $Y_i$ of S&P 500 Index during the period 3 January 1928 to 30 August 1991, and concluded that: (i) The S&P 500 returns $Y_i$ are not i.i.d.; and (ii) $Y_i$ do have some short memory and there is a small amount of predictability in stock returns. See also Verhoeven et al. [33], Caporale and Gil-Alana [6], and Awartani and Corradi [2] for more discussions on the serial dependence among the S&P 500 returns.

Here we shall model this dataset using model (1.1) with $(X_i, Y_i) = (Y_{i-1}, Y_i)$, $i = 1, 2, \ldots, 14372$. Since 13562 (94.4%) out of the 14373 $Y$'s lie within the range $[-0.017, 0.017]$, we only keep those pairs $(X_i, Y_i)$ for which $X_i \in \mathcal{T} := [-0.017, 0.017]$.

We now address the bandwidth selection issue. Popular bandwidth selection methods include plug-in method, cross-validation approach and Fan and Gijbels's [14] residual squares criterion (RSC) method. As demonstrated by the latter paper, the RSC method is efficient in choosing bandwidth in local polynomial regression. Here we shall briefly illustrate the idea. Suppose that we are interested in local linear estimation of $m(x) = \mathbb{E}(Y|X = x)$ based on data $(X_i, Y_i)_{i=1}^n$ drawn from the distribution $(X, Y)$. Let $\mathbf{X}$ be the design matrix whose $i$th row is $(1, X_i - x)$, $1 \le i \le n$, $\mathbf{Y} = (Y_1, Y_2, \ldots, Y_n)^T$ the response vector, and $\mathbf{K} = \mathrm{diag}\{K_{b_n}(X_1 - x), K_{b_n}(X_2 - x), \ldots, K_{b_n}(X_n - x)\}$ the diagonal weighting matrix. The RSC is defined as

$$\mathrm{RSC}(x; b_n) = \frac{1 + 2\lambda}{\mathrm{trace}\{\mathbf{K} - \mathbf{KXS}_n^{-1}\mathbf{X}^T\mathbf{K}\}}(\mathbf{Y} - \mathbf{X}\hat{\beta})^T\mathbf{K}(\mathbf{Y} - \mathbf{X}\hat{\beta}),$$

TABLE 2
*Simulated coverage probabilities for Model 2*

| $h_n$ | 0.16 | 0.18 | 0.20 | 0.22 | 0.24 | 0.26 | 0.28 | 0.30 |
|---|---|---|---|---|---|---|---|---|
| use $q_\alpha$ | 0.9625 | 0.9646 | 0.9653 | 0.9688 | 0.9702 | 0.9721 | 0.9679 | 0.9654 |
| use $q_\alpha^*$ | 0.9435 | 0.9443 | 0.9490 | 0.9534 | 0.9529 | 0.9572 | 0.9525 | 0.9498 |

where $\hat{\beta} = \mathbf{S}_n^{-1} \mathbf{X}^T \mathbf{K} \mathbf{Y}, \mathbf{S}_n = \mathbf{X}^T \mathbf{K} \mathbf{X}, \mathbf{S}_n^* = \mathbf{X}^T \mathbf{K}^2 \mathbf{X}$ and $\lambda$ is the first diagonal element of the matrix $\mathbf{S}_n^{-1} \mathbf{S}_n^* \mathbf{S}_n^{-1}$. The idea is to find $b_n$ to minimize the integrated version of the RSC: $\mathrm{IRSC}(b_n) = \int_{\mathcal{T}} \mathrm{RSC}(x; b_n)\, dx$; see Fan and Gijbels's paper for more details. With their method, we outline our procedure:

(1) Define $\mathrm{IRSC}_\mu(b_n)$ for bandwidth selection in estimating $\mu$.

(2) Obtain $b_n^*$ by minimizing $\mathrm{IRSC}_\mu(b_n)$.

(3) Get local linear regression estimate $\hat{\mu}_{b_n^*}$ and $\hat{\mu}_{\sqrt{2}b_n^*}$.

(4) Apply bias-correction procedure and obtain the final estimate $\hat{\mu}_{b_n^*}^* = 2\hat{\mu}_{b_n^*} - \hat{\mu}_{\sqrt{2}b_n^*}$.

(5) Compute residuals, and define $\mathrm{IRSC}_\sigma(h_n)$ for bandwidth selection in estimating $\sigma^2$.

(6) Obtain $h_n^*$ by minimizing $\mathrm{IRSC}_\sigma(b_n)$.

(7) As in steps (3) and (4), calculate $\hat{\sigma}_{h_n^*}^{2*}$.

We find $b_n^* = 0.013$ with $\mathrm{IRSC}_\mu(b_n^*) = 2.292 \times 10^{-6}$. Furthermore, we compute the relative IRSC: $\mathrm{IRSC}_\mu(b_n)/\mathrm{IRSC}_\mu(b_n^*) = 1.226, 1.194, 1.154, 1.114, 1.078, 1.050$ for $b_n = 0.001, 0.002, \ldots, 0.006$, respectively, and the IRSC curve tends to be flat after $b_n = 0.005$. To avoid over-smoothing and yet control IRSC, we choose $b_n = 0.005$ as our final bandwidth. Similarly, we choose $h_n = 0.006$ with the relative IRSC being $1.063$. We use the cutoff value $q_\alpha^*$ and the grid points $\mathcal{T}_{(30)} = \{-0.017 + 0.034j/29 : j = 0, 1, \ldots, 29\}$ described in Section 5 to construct SCB for $\mu(\cdot)$ and $\sigma^2(\cdot)$ in model (1.1). The estimated $\hat{\nu}_\varepsilon = 6.61$. We have also tried various choices of bandwidths and obtained very similar results.

Interestingly, the 95% SCB for $\mu$ and $\sigma^2$ in Figure 1 suggest that we can accept the two null hypotheses that the regression function $\mu(\cdot)$ is linear and that the volatility function $\sigma^2(\cdot)$ is quadratic. The fitted linear equation is $\hat{\mu}_{\mathrm{linear}}(x) = 0.00022 + 0.138x$ and the fitted quadratic curve is $\hat{\sigma}_{\mathrm{quadratic}}^2(x) = 0.000058 - 0.0011x + 0.257x^2$. Moreover, since the fitted constant line (solid) $\hat{\sigma}_{\mathrm{constant}}^2(x) = 0.000068$ in Figure 1 for $\sigma^2(\cdot)$ is not entirely contained within the 95% SCB, we claim the existence of conditional heteroscedasticity. Finally, we conclude that the following AR(1)–ARCH(1) model is an adequate fit for the S&P 500 Index data:

$$(6.1) \quad Y_i = 0.00022 + 0.138Y_{i-1} + \varepsilon_i \sqrt{0.000058 - 0.0011Y_{i-1} + 0.257Y_{i-1}^2}.$$

The volatility function in (6.1) is slightly different from that of Engle's [10] ARCH model $Y_i = (a^2 + b^2 Y_{i-1}^2)^{1/2} \varepsilon_i$. The latter one implies symmetric conditional volatility while the former one allows for asymmetric volatility: The negative coefficient $-0.0011$ of $Y_{i-1}$ in (6.1) is in line with the empirical observation that bad news (negative return) causes bigger volatility than good news (positive return). The linear term in Model (6.1) suggests some amount of predictability in returns. While it is widely accepted that financial market is efficient and consequently stock prices are not predictable, there is still an increasing number of works that question such a theory. See Ding, Granger and Engle [8], Fama and French [12], Gropp [20] and references therein. On the other hand, since the volatility term

$$\sqrt{0.000058 - 0.0011Y_{i-1} + 0.257Y_{i-1}^2} > 0.5|Y_{i-1}|$$

dominates the linear term, the predictability is fairly weak relative to the noise level.

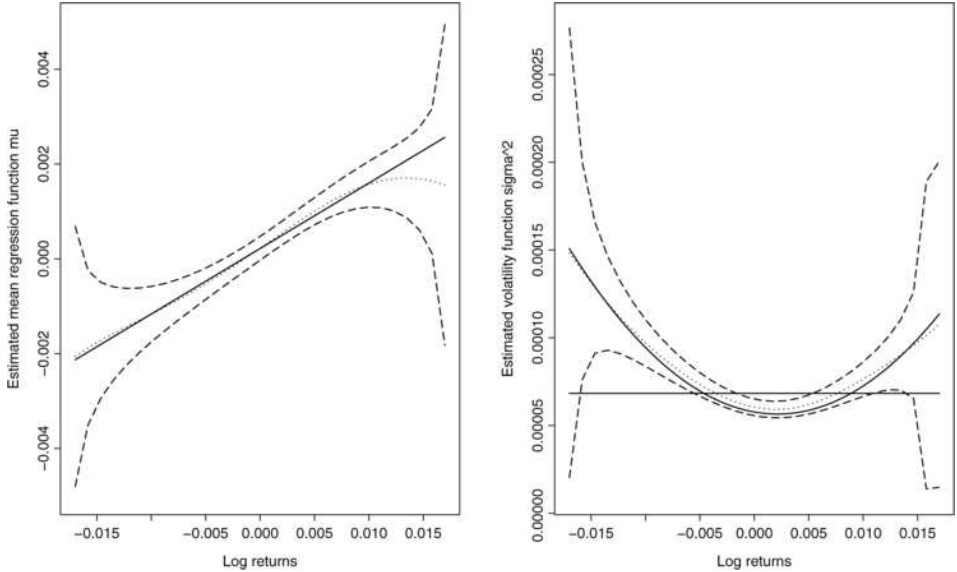

Fig. 1.   *Left: The 95% SCB for regression function $\mu$; dotted, long-dashed and solid lines are the estimated curve $\hat{\mu}^*_{b_n}$, SCB and the fitted linear line $\hat{\mu}_{\text{linear}}(x) = 0.00022 + 0.1382x$, respectively. Right: The 95% SCB for volatility function $\sigma^2$; dotted, long-dashed, solid quadratic and solid constant lines are the estimated curve $\hat{\sigma}^{2*}_{b_n}$, SCB, the fitted quadratic curve $\hat{\sigma}^2_{\text{quadratic}}(x) = 0.000058 - 0.0011x + 0.2567x^2$, and the constant fit $\hat{\sigma}^2_{\text{constant}}(x) = 0.000068$, respectively.*

**7. Proofs.**   This section provides proofs of results stated in Section 3. Recall $\mathcal{F}_i = (\ldots, \eta_{i-1}, \eta_i)$. Let $\mathcal{G}_i = (\ldots, \eta_i, \eta_{i+1}; \varepsilon_i, \varepsilon_{i-1}, \ldots)$. By the assumption in Section 2, $\varepsilon_i$ is independent of $\mathcal{G}_{i-1}$. In the sequel, with a slight abuse of notation we refer $\mathcal{F}_i$ (resp. $\mathcal{G}_i$) as the sigma field generated by $\mathcal{F}_i$ (resp. $\mathcal{G}_i$). Recall (2.3) for $\Xi_n$. Throughout the proofs we assume WLOG that the kernel $K$ has bounded support $[-1, 1]$.

Recall that $f_X(x|\mathcal{F}_{i-1})$ is the conditional density of $X_i$ at $x$ given $\mathcal{F}_{i-1}$. Define

$$(7.1) \qquad I_n(x) = \sum_{i=1}^{n} [f_X(x|\mathcal{F}_{i-1}) - \mathbb{E}f_X(x|\mathcal{F}_{i-1})], \qquad x \in \mathbb{R}.$$

LEMMA 1. *Let $T > 0$ be fixed. Then $\| \sup_{|x| \le T} |I_n(x)| \| = O(\Xi_n^{1/2})$.*

PROOF. By Theorem 1 in Wu [36], $\sup_{x \in \mathbb{R}} [\|I_n(x)\| + \|I'_n(x)\|] = O(\Xi_n^{1/2})$. Then Lemma 1 easily follows since $\sup_{|x| \le T} |I_n(x) - I_n(-T)| \le \int_{-T}^{T} |\partial I_n(x)/\partial x| \, dx$. □

### 7.1. A CLT and a maximal deviation result.

Let $g$ and $h$ be measurable functions such that $h(\varepsilon_0) \in \mathcal{L}^2$ and $\vartheta_h^2 = \mathrm{Var}[h(\varepsilon_0)] > 0$. For $K \in \mathcal{K}$ define

$$
\begin{aligned}
(7.2) \qquad & S_n(x) = \sum_{i=1}^{n} \xi_i(x), \\
& \text{where } \xi_i(x) = \frac{g(X_i)[h(\varepsilon_i) - \mathbb{E}(h(\varepsilon_i))]K_{b_n}(x - X_i)}{\vartheta_h g(x)\sqrt{nb_n}\varphi_K f_X(x)}.
\end{aligned}
$$

In Theorems 4 and 5 below, we shall establish a central limit theorem and a maximal deviation result for $S_n(x)$, respectively. These results are of independent interest and they are essential to the proofs of our main results in Section 3.

THEOREM 4. *Let $x \in \mathbb{R}$ be fixed, $K \in \mathcal{K}$ and $h(\varepsilon_0) \in \mathcal{L}^2$. Assume that $f_X(x) > 0$, $g(x) \ne 0$, and $f_X, g \in \mathcal{C}^0(\{x\}^\epsilon)$ for some $\epsilon > 0$. Further assume that $b_n \to 0$, $nb_n \to \infty$ and $\Xi_n/n^2 \to 0$. Then $S_n(x) \Rightarrow N(0,1)$.*

PROOF. Since $\varepsilon_i$ is independent of $\mathcal{G}_{i-1}$, $\{\xi_i(x)\}_{i=1}^{n}$ form martingale differences with respect to $\mathcal{G}_i$. By the martingale central limit theorem, it suffices to verify the convergence of conditional variance and the Lindeberg condition. Let $\gamma_i = g^2(X_i)K_{b_n}^2(x - X_i)$, $u_i = \gamma_i - \mathbb{E}(\gamma_i|\mathcal{F}_{i-1})$ and $v_i = \mathbb{E}(\gamma_i|\mathcal{F}_{i-1}) - \mathbb{E}(\gamma_i)$. Write

$$(7.3) \qquad \sum_{i=1}^{n} [\gamma_i - \mathbb{E}(\gamma_i)] = M_n + R_n, \qquad \text{where } M_n = \sum_{i=1}^{n} u_i \text{ and } R_n = \sum_{i=1}^{n} v_i.$$

Hereafter we shall call (7.3) the $M/R$-decomposition. Since $\{u_i\}_{i=1}^{n}$ are martingale differences with respect to $\mathcal{F}_i$ and $\mathbb{E}(u_i^2) = O(b_n)$, we have $M_n = O_p(\sqrt{nb_n})$. Recall (7.1) for $I_n(x)$. Since $K$ is bounded and has bounded support and $g \in \mathcal{C}^0(\{x\}^\epsilon)$, by Lemma 1,

$$
\begin{aligned}
\|R_n\| &= b_n \left\| \int_{\mathbb{R}} K^2(u) g^2(x - ub_n) I_n(x - ub_n) \, du \right\| \\
&\leq b_n \int_{\mathbb{R}} K^2(u) g^2(x - ub_n) \|I_n(x - ub_n)\| \, du = O(\Xi_n^{1/2} b_n).
\end{aligned}
\tag{7.4}
$$

Since $b_n \to 0$, $nb_n \to \infty$ and $\Xi_n/n^2 \to 0$, simple calculations show that

$$
\sum_{i=1}^n \mathbb{E}[\xi_i^2(x) | \mathcal{G}_{i-1}] = \frac{M_n + R_n}{nb_n \varphi_K f_X(x) g^2(x)} + \frac{\sum_{i=1}^n \mathbb{E}(\gamma_i)}{nb_n \varphi_K f_X(x) g^2(x)} \xrightarrow{\mathrm{p}} 1.
\tag{7.5}
$$

There exists $c$ such that $\sup_u |g(u) K_{b_n}(x - u)| \leq c$ holds for all sufficiently large $n$. Let $\lambda = \vartheta_h g(x) [\varphi_K f_X(x)]^{1/2}$ and $\bar{h}(\varepsilon_0) = h(\varepsilon_0) - \mathbb{E}[h(\varepsilon_0)]$. For any $s > 0$, by the independence of $X_0$ and $\varepsilon_0$,

$$
\begin{aligned}
&\sum_{i=1}^n \mathbb{E}[\xi_i^2(x) \mathbf{1}_{|\xi_i(x)| \geq s}] \\
&\quad = \frac{1}{\lambda^2 b_n} \mathbb{E}[g^2(X_0) K_{b_n}^2(x - X_0) \bar{h}^2(\varepsilon_0) \mathbf{1}_{|g(X_0) K_{b_n}(x - X_0) \bar{h}(\varepsilon_0)| \geq \lambda s \sqrt{nb_n}}] \\
&\quad \leq \frac{1}{\lambda^2 b_n} \mathbb{E}[g^2(X_0) K_{b_n}^2(x - X_0) \bar{h}^2(\varepsilon_0) \mathbf{1}_{|\bar{h}(\varepsilon_0)| \geq \lambda s \sqrt{nb_n}/c}] \\
&\quad = \frac{1}{\lambda^2 b_n} \mathbb{E}[g^2(X_0) K_{b_n}^2(x - X_0)] \times \mathbb{E}[\bar{h}^2(\varepsilon_0) \mathbf{1}_{|\bar{h}(\varepsilon_0)| \geq \lambda s \sqrt{nb_n}/c}] \to 0
\end{aligned}
$$

in view of $nb_n \to \infty$ and $\bar{h}(\varepsilon_0) \in \mathcal{L}^2$. So the Lindeberg condition holds. $\square$

Recall Theorem 2 for the definitions of $m_n$ and $\mathcal{T}_n$. Let $S_n(x)$ be as in (7.2). Theorem 5 below provides a maximal deviation result for $\sup_{x \in \mathcal{T}_n} |S_n(x)|$. Results of this type are essential to the construction of SCB (cf. Bickel and Rosenblatt [4], Johnston [26] and Eubank and Speckman [11] among others). To obtain a maximal deviation result under dependence, we shall apply Grama and Haeusler's [19] martingale moderate deviation theorem.

THEOREM 5. *Let* $K \in \mathcal{K}$ *and* $h(\varepsilon_0) \in \mathcal{L}^3$. *Assume* $\inf_{u \in \mathcal{T}} f_X(u) > 0$, $g(x) \neq 0, x \in \mathcal{T}$, *and* $f_X, g \in \mathcal{C}^4(\mathcal{T}^\epsilon)$ *for some* $\epsilon > 0$. *Further assume that*

$$
b_n^{4/3} \log n + \frac{(\log n)^3}{nb_n^3} + \frac{\Xi_n (\log n)^2}{n^2 b_n^{4/3}} \to 0.
\tag{7.6}
$$

*Let* $B_n(z)$ *and* $\mathcal{T}_n$ *be as in Theorem 2. Then*

$$
\lim_{n \to \infty} \mathbb{P}\left\{ \sup_{x \in \mathcal{T}_n} |S_n(x)| \leq B_{m_n}(z) \right\} = e^{-2e^{-z}}, \qquad z \in \mathbb{R}.
\tag{7.7}
$$

PROOF. Recall (7.2) for $\xi_i(x)$. For fixed $k \in \mathbb{N}$ and mutually different integers $0 \le j_1, j_2, \ldots, j_k \le m_n$, let the $k$-dimensional vector $\zeta_i = [\xi_i(x_{j_1}), \xi_i(x_{j_2}), \ldots, \xi_i(x_{j_k})]^T$ and $S_{n,k} = \sum_{i=1}^n \zeta_i = [S_n(x_{j_1}), S_n(x_{j_2}), \ldots, S_n(x_{j_k})]^T$. Here $^T$ denotes transpose. Then $\{\zeta_i\}_{i=1}^n$ are $k$-dimensional martingale differences with respect to $\mathcal{G}_i$. Let $Q$ denote the quadratic characteristic matrix of $S_{n,k}$, namely,

$$(7.8) \qquad Q = \sum_{i=1}^n \mathbb{E}(\zeta_i \zeta_i^T | \mathcal{G}_{i-1}) := (Q_{rr'})_{1 \le r, r' \le k}.$$

Let $\tau_{rr'} = \varphi_K g(x_{j_r}) g(x_{j_{r'}}) [f_X(x_{j_r}) f_X(x_{j_{r'}})]^{1/2}$ and write

$$Q_{rr'} = \sum_{i=1}^n \mathbb{E}[\xi_i(x_{j_r}) \xi_i(x_{j_{r'}}) | \mathcal{G}_{i-1}]$$

$$= \frac{1}{n b_n \tau_{rr'}} \sum_{i=1}^n g^2(X_i) K_{b_n}(x_{j_r} - X_i) K_{b_n}(x_{j_{r'}} - X_i).$$

For $r \ne r'$, since $|x_{j_r} - x_{j_{r'}}| \ge 2b_n$ and $K$ has support $[-1, 1]$, $Q_{rr'} = 0$. For $r = r'$, we use the $M/R$-decomposition technique in (7.3). Define

$$\alpha_i(r) = g^2(X_i) K_{b_n}^2(x_{j_r} - X_i) - \mathbb{E}[g^2(X_i) K_{b_n}^2(x_{j_r} - X_i) | \mathcal{F}_{i-1}],$$

$$\beta_i(r) = \mathbb{E}[g^2(X_i) K_{b_n}^2(x_{j_r} - X_i) | \mathcal{F}_{i-1}] - \mathbb{E}[g^2(X_i) K_{b_n}^2(x_{j_r} - X_i)].$$

Since $\{\alpha_i(r)\}_{i=1}^n$ form martingale differences with respect to $\mathcal{F}_i$, we have

$$(7.9) \qquad \left\| \sum_{i=1}^n \alpha_i(r) \right\|_2 = \left[ \sum_{i=1}^n \| \alpha_i(r) \|_2^2 \right]^{1/2} = O(\sqrt{nb_n}),$$

uniformly over $r$. By Schwarz's inequality and Lemma 1, as in (7.4), we have

$$\left\| \sum_{i=1}^n \beta_i(r) \right\|_2$$

$$(7.10) \qquad = b_n \left\| \int_{\mathbb{R}} K^2(u) g^2(x_{j_r} - ub_n) I_n(x_{j_r} - ub_n) \, du \right\|_2$$

$$\le b_n \int_{\mathbb{R}} K^2(u) g^2(x_{j_r} - ub_n) \| I_n(x_{j_r} - ub_n) \|_2 \, du = O(\Xi_n^{1/2} b_n),$$

uniformly over $r$. Since $f_X, g \in \mathcal{C}^4(\mathcal{T}^\epsilon)$ and $K \in \mathcal{K}$, by Taylor's expansion,

$$(7.11) \qquad \left| \sum_{i=1}^n \mathbb{E}[g^2(X_i) K_{b_n}^2(x_{j_r} - X_i)] - n b_n \tau_{rr} \right| = O(n b_n^3).$$

Let $\delta_n = (nb_n)^{-1/2} + b_n^4 + \Xi_n^{1/2}/n$. By (7.9), (7.10) and (7.11),

$$
\begin{aligned}
\|Q_{rr} - 1\|_{3/2} \leq \frac{1}{nb_n\tau_{rr}} & \left\{ \left\| \sum_{i=1}^n \alpha_i(r) \right\|_{3/2} + \left\| \sum_{i=1}^n \beta_i(r) \right\|_{3/2} \right. \\
& \left. + \left\| \sum_{i=1}^n \mathbb{E}[g^2(X_i)K_{b_n}^2(x_{j_r} - X_i)] - nb_n\tau_{rr} \right\|_{3/2} \right\} \\
& = O(\delta_n)
\end{aligned}
$$
(7.12)

uniformly over $r$. Let $\mathbb{I}_k = \mathrm{diag}(1,1,\ldots,1) = (u_{rr'})_{1 \leq r,r' \leq k}$ be the $k \times k$ identity matrix. Then $\mathbb{E}|Q_{rr'} - u_{rr'}|^{3/2} = O(\delta_n^{3/2})$, uniformly over $1 \leq r, r' \leq k$. It is easily seen that $\sum_{i=1}^n \mathbb{E}|\xi_i(x_{j_r})|^3 = O[(nb_n)^{-1/2}]$ uniformly over $1 \leq r \leq k$. Then $\sum_{i=1}^n \mathbb{E}|\xi_i(x_{j_r})|^3 + \mathbb{E}|Q_{rr'} - u_{rr'}|^{3/2} = O(\Lambda_n)$, where $\Lambda_n = (nb_n)^{-1/2} + b_n^3 + \Xi_n^{3/4}/n^{3/2}$.

Under (7.6), elementary calculations show that $[1 + B_{m_n}(z)]^4 \exp[B_{m_n}^2(z)/2] \times \Lambda_n \to 0$ for fixed $z$. Let $A_j$ denote the event $\{|S_n(x_j)| > B_{m_n}(z)\}$ and $E_{m_n} = \bigcup_{j=0}^{m_n} A_j$; let $N_1, N_2, \ldots,$ be i.i.d. stand normals. By Theorem 1 in Grama and Haeusler [19],

$$
\mathbb{P}\left[ \bigcap_{r=1}^k A_{j_r} \right] = \mathbb{P}\left[ \bigcap_{r=1}^k \{|N_r| > B_{m_n}(z)\} \right][1 + o(1)] = \left( \frac{2e^{-z}}{m_n} \right)^k [1 + o(1)],
$$
(7.13)

in view of $\mathbb{P}(N_1 > x) = [1 + o(1)]\phi(x)/x$ as $x \to \infty$, where $\phi$ is the standard normal density function. Notice that $\mathbb{P}\{\sup_{x \in \mathcal{T}_n} |S_n(x)| > B_{m_n}(z)\} = \mathbb{P}(E_{m_n})$. By the inclusion-exclusion inequality, we have, for large enough $n$,

$$
\begin{aligned}
\mathbb{P}(E_{m_n}) & \leq \sum_{j=0}^{m_n} \mathbb{P}[A_j] - \sum_{j_1 < j_2 \leq m_n} \mathbb{P}[A_{j_1} \cap A_{j_2}] + \cdots \\
& \quad + \sum_{j_1 < j_2 < \cdots < j_{2k-1} \leq m_n} \mathbb{P}\left[ \bigcap_{r=1}^{2k-1} A_{j_r} \right] \\
& = \sum_{r=1}^{2k-1} (-1)^{r-1} \binom{m_n + 1}{r} \left( \frac{2e^{-z}}{m_n} \right)^r [1 + o(1)] \\
& = -\sum_{r=1}^{2k-1} \frac{(-2e^{-z})^r}{r!} [1 + o(1)].
\end{aligned}
$$
(7.14)

Thus, let $\tau_j = \sum_{r=1}^j [-2\exp(-z)]^r/r!$, we have $\limsup_{n \to \infty} \mathbb{P}(E_{m_n}) \leq -\tau_{2k-1}$. Similarly, $\liminf_{n \to \infty} \mathbb{P}(E_{m_n}) \geq -\tau_{2k}$. Since $\lim_{k \to \infty} \tau_k = 1 - \exp[-2\exp(-z)]$, (7.7) follows. □

7.2. *Proof of the results in Section* 3. Before proving Theorems 1 and 2, we first note that, by the definition of $\hat{\mu}_{b_n}(x)$ and $\hat{f}_X(x)$ in (3.1),

$$(7.15) \quad \frac{\sqrt{nb_n}}{\sigma(x)\sqrt{\varphi_K}}\sqrt{\hat{f}_X(x)}\{\hat{\mu}_{b_n}(x) - \mu(x) - \omega_n(x)U_n(x)\} = V_n(x)\sqrt{\omega_n(x)},$$

where $\omega_n(x) = f_X(x)/\hat{f}_X(x)$,

$$(7.16) \quad U_n(x) = \frac{1}{nb_nf_X(x)}\sum_{i=1}^{n}K_{b_n}(x - X_i)[\mu(X_i) - \mu(x)],$$

$$(7.17) \quad V_n(x) = \frac{1}{\sigma(x)\sqrt{nb_n\varphi_Kf_X(x)}}\sum_{i=1}^{n}\sigma(X_i)\varepsilon_iK_{b_n}(x - X_i).$$

In (7.15), we can view $U_n(x)$ [resp. $V_n(x)$] as the bias (resp. stochastic) part of $\hat{\mu}_{b_n}(x) - \mu(x)$. The stochastic part $V_n(x)$ is treated in Proposition 4 and Theorem 5. The following Lemma 2 concerns $\omega_n(x)$ and the bias part $U_n(x)$.

LEMMA 2. *Let $K \in \mathcal{K}$.*

(i) *Recall the definition of $\hat{f}_X(x)$ in (3.1). Assume that $f_X \in \mathcal{C}^4(\mathcal{T}^\epsilon)$ for some $\epsilon > 0$. Then*

$$(7.18) \quad \sup_{x\in\mathcal{T}}|\hat{f}_X(x) - f_X(x)| = O_p(q_n),$$
$$where \ q_n = \sqrt{\log n/(nb_n)} + b_n^2 + \Xi_n^{1/2}/n.$$

(ii) *Recall the definition of $U_n(x)$ in (7.16). Let $\rho_\mu(x)$ be as in Theorem 1. Assume that $f_X, \mu \in \mathcal{C}^4(\mathcal{T}^\epsilon)$ for some $\epsilon > 0$ and $\inf_{x\in\mathcal{T}}f_X(x) > 0$. Then*

$$(7.19) \quad \sup_{x\in\mathcal{T}}|U_n(x) - b_n^2\psi_K\rho_\mu(x)| = O_p(r_n),$$
$$where \quad r_n = \sqrt{b_n\log n/n} + b_n^4 + \Xi_n^{1/2}b_n/n.$$

(iii) *Let $g, h$ be measurable functions such that $h(\varepsilon_0) \in \mathcal{L}^q$ for some $q \geq 2$ and $g \in \mathcal{C}^0(\mathcal{T}^\epsilon)$ for some $\epsilon > 0$. Then*

$$(7.20) \quad \sup_{x\in\mathcal{T}}\frac{1}{nh_n}\left|\sum_{i=1}^{n}K_{b_n}(x - X_i)g(X_i)[h(\varepsilon_i) - \mathbb{E}h(\varepsilon_i)]\right| = O_p(\chi_n(q)),$$
$$where \ \chi_n(q) = \sqrt{\log n/(nb_n)} + n^{-q/4}b_n^{-q/4-1}(\log n)^{q/4-1/2}.$$

PROOF. (i) We shall use the $M/R$-decomposition technique in (7.3). By the chain argument in Lemma 4 in Zhao and Wu [39], we can show that

$$(7.21) \quad \sup_{x\in\mathcal{T}}\left|\sum_{i=1}^{n}\{K_{b_n}(x - X_i) - \mathbb{E}[K_{b_n}(x - X_i)|\mathcal{F}_{i-1}]\}\right| = O_p(\sqrt{nb_n\log n}).$$

By Lemma 1, since $K$ has bounded support, we have

$$
\begin{aligned}
(7.22) \quad & \sum_{i=1}^{n} \{ \mathbb{E}[K_{b_n}(x - X_i) | \mathcal{F}_{i-1}] - \mathbb{E}[K_{b_n}(x - X_i)] \} \\
& = b_n \int_{\mathbb{R}} K(u) f_X(x - u b_n) I_n(x - u b_n) \, du = O_{\mathrm{p}}(\Xi_n^{1/2} b_n).
\end{aligned}
$$

Since $K \in \mathcal{K}$, by Taylor's expansion, $\sum_{i=1}^{n} \mathbb{E}[K_{b_n}(x - X_i)] - n b_n f_X(x) = O(n b_n^3)$. Therefore, (i) follows from (7.21) and (7.22).

(ii) Let $\gamma_i(x) = K_{b_n}(x - X_i)[\mu(X_i) - \mu(x)]$. Again we use the $M/R$-decomposition technique in (7.3). By the chain argument in Lemma 4 in Zhao and Wu [39], we can show that $\sup_{x \in \mathcal{T}} | \sum_{i=1}^{n} \{ \gamma_i(x) - \mathbb{E}[\gamma_i(x) | \mathcal{F}_{i-1}] \} | = O_{\mathrm{p}}[(n b_n^3 \log n)^{-1/2}]$. As in (7.22), we can show that $\sup_{x \in \mathcal{T}} | \sum_{i=1}^{n} \{ \mathbb{E}[\gamma_i(x) | \mathcal{F}_{i-1}] - \mathbb{E}[\gamma_i(x)] \} | = O(\Xi_n^{1/2} b_n^2)$. Since $K \in \mathcal{K}$, elementary calculations show that $\sum_{i=1}^{n} \mathbb{E}[\gamma_i(x)] = n b_n^3 f_X(x) \psi_K \rho_\mu(x) + O(n b_n^5)$. So (ii) follows.

(iii) We shall only consider the special case of $h(u) = u$ since other cases follow similarly. Let $c_n = (n b_n / \log n)^{1/2}$ and define

$$
\overline{D}_n(x) = \sum_{i=1}^{n} \overline{d}_i(x),
$$

$$
\text{where } \overline{d}_i(x) = K_{b_n}(x - X_i) g(X_i)[\varepsilon_i \mathbf{1}_{|\varepsilon_i| > c_n} - \mathbb{E}(\varepsilon_i \mathbf{1}_{|\varepsilon_i| > c_n})],
$$

$$
\underline{D}_n(x) = \sum_{i=1}^{n} \underline{d}_i(x),
$$

$$
\text{where } \underline{d}_i(x) = K_{b_n}(x - X_i) g(X_i)[\varepsilon_i \mathbf{1}_{|\varepsilon_i| \le c_n} - \mathbb{E}(\varepsilon_i \mathbf{1}_{|\varepsilon_i| \le c_n})]/c_n.
$$

Note that, for each fixed $x$, $\{ \overline{d}_i(x) \}_{i=1}^{n}$ and $\{ \underline{d}_i(x) \}_{i=1}^{n}$ form martingale differences with respect to $\mathcal{G}_i$, and $\mathbb{E}(\varepsilon_i^2 \mathbf{1}_{|\varepsilon_i| > c_n}) \le \mathbb{E}(|\varepsilon_i|^q)/c_n^{q-2} = O(c_n^{2-q})$ for $q \ge 2$. Simple calculations show that $\| \overline{D}_n(x) \|^2 = \sum_{i=1}^{n} \| \overline{d}_i(x) \|^2 = O(n b_n c_n^{2-q})$ and $\| \partial \overline{D}_n(x) / \partial x \|^2 = O(n c_n^{2-q} / b_n)$, uniformly over $x \in \mathcal{T}$. Since $\sup_{x \in \mathcal{T}} | \overline{D}_n(x) - \overline{D}_n(T_1) | \le \int_{\mathcal{T}} | \partial \overline{D}_n(u) / \partial u | \, du$, by Schwarz's inequality, we have $\mathbb{E}[\sup_{x \in \mathcal{T}} | \overline{D}_n(x) |^2] = O(n c_n^{2-q} / b_n)$. Since $\{ \underline{d}_i(x) \}_{i=1}^{n}$ are uniformly bounded martingale differences, by the argument in the proof of Lemma 4 in Zhao and Wu [39], $\sup_{x \in \mathcal{T}} | \underline{D}_n(x) | = O_{\mathrm{p}}[(n b_n \log n)^{1/2} / c_n]$. Note that $\mathbb{E}(\varepsilon_i) = 0$. Then $| \sum_{i=1}^{n} K_{b_n}(x - X_i) g(X_i) \varepsilon_i | \le | \overline{D}_n(x) | + c_n | \underline{D}_n(x) |$. So (iii) follows.  $\square$

PROOF OF THEOREMS 1 AND 2.  Recall the definition of $\omega_n(x)$, $U_n(x)$ and $V_n(x)$ in (7.15). Applying the $M/R$-decomposition technique in (7.3), we can show that, for fixed $x$, $\hat{f}_X(x) - f_X(x) = O_{\mathrm{p}}[(n b_n)^{-1/2} + b_n^2 + \Xi_n^{1/2} / n]$ and $U_n(x) - b_n^2 \psi_K \rho_\mu(x) = O_{\mathrm{p}}[(b_n / n)^{1/2} + b_n^2 + \Xi_n^{1/2} b_n / n]$. Then $\omega_n(x) = 1 + O_{\mathrm{p}}[(n b_n)^{-1/2} + b_n^2 + \Xi_n^{1/2} / n]$. By Theorem 4, $V_n(x) \Rightarrow N(0, 1)$. Under condition (3.2), Theorem 1 then follows from Slutsky's theorem.

By Lemma 2, $\sup_{x \in \mathcal{T}} |\omega_n(x) - 1| = O_p(q_n)$ and $\sup_{x \in \mathcal{T}} |U_n(x) - b_n^2 \psi_K \times \rho_\mu(x)| = O_p(r_n)$, where $q_n$ and $r_n$ are as in (7.18) and (7.19), respectively. Under condition (3.3), simple calculations show that $q_n \log n + (nb_n \log n)^{1/2} r_n \to 0$ and that (7.6) holds. So Theorem 2 follows from the decomposition (7.15) in view of Slutsky's theorem and Theorem 5. $\square$

PROOF OF PROPOSITION 1 AND THEOREM 3.    Let $q_n, r_n$ and $\chi_n(q)$ be as in (7.18), (7.19) and (7.20) in Lemma 2, respectively. Accordingly, we define $\tilde{q}_n, \tilde{r}_n$ and $\tilde{\chi}_n(q)$ with $b_n$ in $q_n, r_n$ and $\chi_n(q)$ replaced by $h_n$. For example, $\tilde{q}_n = [\log n/(nh_n)]^{1/2} + h_n^2 + \Xi_n^{1/2}/n$. Recall the definition of $\omega_n(x)$ and $U_n(x)$ in (7.15). Define

$$W_n(x) = \sum_{i=1}^n \frac{\sigma(X_i)\varepsilon_i K_{b_n}(x - X_i)}{nb_n f_X(x)}, \qquad W_n^*(x) = \sum_{i=1}^n \frac{\sigma(X_i)\varepsilon_i K_{b_n}^*(x - X_i)}{nb_n f_X(x)},$$

where $K^*(u) = 2K(u) - K(u/\sqrt{2})/\sqrt{2}$. Applying Lemma 2(iii) with $h(u) = u$ and $q = 6$, we have $\sup_{x \in \mathcal{T}} |W_n^*(x)| = O_p[\chi_n(6)]$. So, by Lemma 2(i) and (ii), elementary calculations show that, uniformly over $x \in \mathcal{T}$,

$$
\begin{aligned}
(7.23) \qquad \hat{\mu}_{b_n}(x) - \mu(x) &= \omega_n(x) U_n(x) + \omega_n(x) W_n(x) \\
&= b_n^2 \psi_K \rho_\mu(x) + W_n(x) + O_p(\Delta_n),
\end{aligned}
$$

where $\Delta_n = r_n + q_n[b_n^2 + \chi_n(6)]$. Consequently, $\hat{\mu}_{b_n}^*(x) = \mu(x) + W_n^*(x) + O_p(\Delta_n)$. Let $\tilde{f}_X$ be as in (3.8). By definition,

$$
\begin{aligned}
(7.24) \qquad \hat{\sigma}_{h_n}^2(x) &= \frac{1}{nh_n \tilde{f}_X(x)} \sum_{i=1}^n [\sigma(X_i)\varepsilon_i - W_n^*(X_i) + O_p(\Delta_n)]^2 K_{h_n}(x - X_i) \\
&= \frac{T_n(x)}{nh_n \tilde{f}_X(x)} + O_p\left[\frac{L_n(x)}{nh_n} + \frac{\Delta_n}{nh_n} J_n(x) + \chi_n^2(6) + \Delta_n^2\right],
\end{aligned}
$$

where

$$T_n(x) = \sum_{i=1}^n \sigma^2(X_i)\varepsilon_i^2 K_{h_n}(x - X_i),$$

$$L_n(x) = \sum_{i=1}^n \sigma(X_i)\varepsilon_i W_n^*(X_i) K_{h_n}(x - X_i)$$

and

$$J_n(x) = \sum_{i=1}^n \sigma(X_i)|\varepsilon_i| K_{h_n}(x - X_i).$$

By the argument in Lemma 2(i), we can show that $\sum_{i=1}^{n} \sigma(X_i)K_{h_n}(x - X_i) = O_p[nh_n(1 + \tilde{q}_n)]$, uniformly over $x \in \mathcal{T}$. By Lemma 2(iii) with $h(u) = |u|$ and $q = 6$, we have uniformly over $x \in \mathcal{T}$ that

$$
\begin{aligned}
(7.25) \quad J_n(x) &= \sum_{i=1}^{n} \sigma(X_i)[|\varepsilon_i| - \mathbb{E}|\varepsilon_i|]K_{h_n}(x - X_i) \\
&\quad + \mathbb{E}|\varepsilon_0| \sum_{i=1}^{n} \sigma(X_i)K_{h_n}(x - X_i) \\
&= O_p\{nh_n[1 + \tilde{q}_n + \tilde{\chi}_n(6)]\}.
\end{aligned}
$$

Since $\varepsilon_i$ is independent of $\mathcal{G}_{i-1}$ and $\mathbb{E}(\varepsilon_i) = 0$, simple calculation shows that

$$
\begin{aligned}
\|L_n(x)\|^2 &= \mathbb{E}\left[\sum_{i,j=1}^{n} \varepsilon_i \varepsilon_j \frac{\sigma(X_i)\sigma(X_j)}{nb_n f_X(X_i)} K_{b_n}^*(X_i - X_j)K_{h_n}(x - X_i)\right]^2 \\
&= O(h_n/b_n^2),
\end{aligned}
$$

uniformly over $x \in \mathcal{T}$. Likewise, $\|\partial L_n(x)/\partial x\|^2 = O[1/(h_n b_n^2)]$. Since $\sup_{x \in \mathcal{T}} |L_n(x) - L_n(T_1)| \leq \int_{\mathcal{T}} |\partial L_n(u)/\partial u| \, du$, by Schwarz's inequality, $\sup_{x \in \mathcal{T}} |L_n(x)| = O_p[1/(h_n^{1/2} b_n)]$. Recall that $\mathbb{E}(\varepsilon_0^2) = 1$. Write

$$
(7.26) \quad \frac{T_n(x)}{nh_n \hat{f}_X(x)} - \sigma^2(x) = \frac{D_n + E_n(x)}{nh_n \tilde{f}_X(x)},
$$

where

$$
\begin{aligned}
D_n(x) &= \sum_{i=1}^{n} [\sigma^2(X_i) - \sigma^2(x)]K_{h_n}(x - X_i), \\
E_n(x) &= \sum_{i=1}^{n} \sigma^2(X_i)[\varepsilon_i^2 - \mathbb{E}(\varepsilon_i^2)]K_{h_n}(x - X_i).
\end{aligned}
$$

Let $\tilde{\omega}_n(x) = f_X(x)/\tilde{f}_X(x)$. By Lemma 2(i), $\sup_{x \in \mathcal{T}} |\tilde{\omega}_n(x) - 1| = O_p(\tilde{q}_n)$. Let $\rho_\sigma(x)$ be as in Theorem 3. As in Lemma 2(ii), we can show that

$$
(7.27) \quad \sup_{x \in \mathcal{T}} \left| \frac{D_n(x)}{nh_n \tilde{f}_X(x)} - h_n^2 \psi_K \rho_\sigma(x) \right| = O_p(\tilde{r}_n).
$$

Thus, by (7.24), (7.25), (7.26) and (7.27), we have

$$
\begin{aligned}
(7.28) \quad \hat{\sigma}_{h_n}^2(x) - \sigma^2(x) - h_n^2 \psi_K \rho_\sigma(x) &= \frac{E_n(x)}{nh_n \tilde{f}_X(x)} + O_p(\ell_n), \\
\text{where } \ell_n &= (nh_n^{3/2} b_n)^{-1} + \tilde{r}_n + [1 + \tilde{q}_n + \tilde{\chi}_n(6)]\Delta_n + \chi_n^2(6) + \Delta_n^2.
\end{aligned}
$$

Applying Lemma 2(iii) with $h(x) = x^2$ and $q = 3$, we have $\sup_{x \in \mathcal{T}} |E_n(x)| = O_p[nh_n \tilde{\chi}_n(3)]$. When $h_n \asymp b_n$ and condition (3.9) is satisfied, Proposition

**1** follows from (7.28) by simplifying $\ell_n + h_n^2 + \tilde{\chi}_n(3)$. Since there are no essential difficulties, we omit the details.

For the proof of Theorem **3**, let $\varsigma_n(x) = \sigma^2(x)/\hat{\sigma}_{h_n}^2(x)$. Since $\mathbb{E}(\varepsilon_0) = 0$, $\mathbb{E}(\varepsilon_0^2) = 1$ and $\varepsilon_0$ has continuous density, we have $\nu_\varepsilon = \mathbb{E}(\varepsilon^4) - 1 > 0$. By (7.28),

$$
\begin{aligned}
(7.29) \quad & \frac{\sqrt{nh_n}}{\sqrt{\varphi_K \nu_\varepsilon}} \frac{[\tilde{f}_X(x)]^{1/2}}{\hat{\sigma}_{h_n}^2(x)} [\hat{\sigma}_{h_n}^2(x) - \sigma^2(x) - h_n^2 \psi_K \rho_\sigma(x)] \\
& = \varsigma_n(x) \sqrt{\tilde{\omega}_n(x)} \frac{E_n(x)}{\sigma^2(x)\sqrt{nh_n \nu_\varepsilon \varphi_K f_X(x)}} + O_p(\sqrt{nh_n}\ell_n).
\end{aligned}
$$

By Proposition **1** and (3.10), $\sup_{x \in \mathcal{T}} |\varsigma_n(x) - 1| = o_p(1/\log n)$. Also, it is easy to check that $\sup_{x \in \mathcal{T}} |\tilde{\omega}_n(x) - 1| = o_p(1/\log n)$ and $(nh_n \log n)^{1/2}\ell_n \to 0$. Thus, Theorem **3** follows from Theorem **5** via Slutsky's theorem. $\square$

PROOF OF PROPOSITION **2**. As shown in the proof of Proposition **1**, we have $\sup_{x \in \mathcal{T}} |\hat{\mu}_{b_n}^* - \mu(x)| = O_p[\Delta_n + \chi_n(6)]$. By (7.28), we have $\sup_{x \in \mathcal{T}} |\hat{\sigma}_{h_n}^{2*}(x) - \sigma^2(x)| = O_p[\tilde{\chi}_n(3) + \ell_n]$. Therefore,

$$
\begin{aligned}
(7.30) \quad \sum_{i=1}^n \hat{\varepsilon}_i \mathbf{1}_{X_i \in \mathcal{T}} & = \sum_{i=1}^n \left[\frac{\sigma(X_i)\varepsilon_i + \mu(X_i) - \hat{\mu}_{b_n}^*(X_i)}{\hat{\sigma}_{h_n}^*(X_i)}\right]^4 \mathbf{1}_{X_i \in \mathcal{T}} \\
& = \sum_{i=1}^n \varepsilon_i^4 \mathbf{1}_{X_i \in \mathcal{T}} + O\{n[\tilde{\chi}_n(3) + \chi_n(6) + \ell_n + \Delta_n]\}.
\end{aligned}
$$

By the independence of $\varepsilon_i$ and $\mathcal{G}_{i-1}$, $\{[\varepsilon_i^4 - \mathbb{E}(\varepsilon_i^4)]\mathbf{1}_{X_i \in \mathcal{T}}\}_{i=1}^n$ form martingale differences with respect to $\mathcal{G}_i$. Since $\varepsilon_i \in \mathcal{L}^6$, $\|\sum_{i=1}^n [\varepsilon_i^4 - \mathbb{E}(\varepsilon_i^4)]\mathbf{1}_{X_i \in \mathcal{T}}\|_{3/2} = O(n^{2/3})$. Furthermore, by applying the $M/R$-decomposition technique in (7.3), we can show that

$$
\begin{aligned}
& \sum_{i=1}^n [\mathbf{1}_{X_i \in \mathcal{T}} - \mathbb{E}(\mathbf{1}_{X_i \in \mathcal{T}})] \\
(7.31) \quad & = \sum_{i=1}^n [\mathbf{1}_{X_i \in \mathcal{T}} - \mathbb{E}(\mathbf{1}_{X_i \in \mathcal{T}} | \mathcal{F}_{i-1})] + \sum_{i=1}^n [\mathbb{E}(\mathbf{1}_{X_i \in \mathcal{T}} | \mathcal{F}_{i-1}) - \mathbb{E}(\mathbf{1}_{X_i \in \mathcal{T}})] \\
& = O_p(\sqrt{n} + \Xi_n^{1/2}).
\end{aligned}
$$

Thus, the desired result follows from (7.30) and (7.31) via elementary manipulations. $\square$

**Acknowledgments.** We are grateful to the referees and an Associate Editor for their many helpful comments. We would also like to thank Erich Haeusler for clarifications on martingale moderate deviations, and Mathias

Drton, Steven Lalley, Peter McCullagh, Mary Sara McPeek, Michael Stein, Matthew Stephens, Mei Wang and Michael Wichura for comments and discussions.

## REFERENCES

[1] Aït-Sahalia, Y. (1996). Nonparametric pricing of interest rate derivative securities. *Econometrica* **64** 527–560.

[2] Awartani, B. M. A. and Corradi, V. (2005). Predicting the volatility of the S&P-500 stock index via GARCH models: The role of asymmetries. *Int. J. Forecasting* **21** 167–183.

[3] Bandi, F. M. and Phillips, P. C. B. (2003). Fully nonparametric estimation of scalar diffusion models. *Econometrica* **71** 241–283. MR1956859

[4] Bickel, P. J. and Rosenblatt, M. (1973). On some global measures of the deviations of density function estimates. *Ann. Statist.* **1** 1071–1095. MR0348906

[5] Bühlmann, P. (1998). Sieve bootstrap for smoothing in nonstationary time series. *Ann. Statist.* **26** 48–83. MR1611804

[6] Caporale, G. M. and Gil-Alana, L. A. (2004). Long range dependence in daily stock returns. *Appl. Financ. Econ.* **14** 375–383.

[7] Cummins, D. J., Filloon, T. G. and Nychka, D. (2001). Confidence intervals for nonparametric curve estimates: Toward more uniform pointwise coverage. *J. Amer. Statist. Assoc.* **96** 233–246. MR1952734

[8] Ding, Z., Granger, C. W. J. and Engle, R. F. (1993). A long memory property of stock market returns and a new model. *J. Empirical Finance* **1** 83–106.

[9] Dümbgen, L. (2003). Optimal confidence bands for shape-restricted curves. *Bernoulli* **9** 423–449. MR1997491

[10] Engle, R. F. (1982). Autoregressive conditional heteroscedasticity with estimates of the variance of U.K. inflation. *Econometrica* **50** 987–1008. MR0666121

[11] Eubank, R. L. and Speckman, P. L. (1993). Confidence bands in nonparametric regression. *J. Amer. Statist. Assoc.* **88** 1287–1301. MR1245362

[12] Fama, E. F. and French, K. R. (1988). Permanent and temporary components of stock prices. *J. Polit. Economy* **96** 246–273.

[13] Fan, J. (2005). A selective overview of nonparametric methods in financial econometrics. *Statist. Sci.* **20** 317–337. MR2210224

[14] Fan, J. and Gijbels, I. (1995). Data-driven bandwidth selection in local polynomial fitting: Variable bandwidth and spatial adaptation. *J. Roy. Statist. Soc. Ser. B* **57** 371–394. MR1323345

[15] Fan, J. and Gijbels, I. (1996). *Local Polynomial Modeling and Its Applications.* Chapman and Hall, London. MR1383587

[16] Fan, J. and Yao, Q. (1998). Efficient estimation of conditional variance functions in stochastic regression. *Biometrika* **85** 645–660. MR1665822

[17] Fan, J. and Yao, Q. (2003). *Nonlinear Time Series*: *Nonparametric and Parametric Methods.* Springer, New York. MR1964455

[18] Fan, J., Zhang, C. and Zhang, J. (2001). Generalized likelihood ratio statistics and Wilks phenomenon. *Ann. Statist.* **29** 153–193. MR1833962

[19] Grama, I. G. and Haeusler, E. (2006). An asymptotic expansion for probabilities of moderate deviations for multivariate martingales. *J. Theoret. Probab.* **19** 1–44. MR2256478

[20] Gropp, J. (2004). Mean reversion of industry stock returns in the U.S., 1926–1998. *J. Empirical Finance* **11** 537–551.

[21] HAGGAN, V. and OZAKI, T. (1981). Modelling nonlinear random vibrations using an amplitude-dependent autoregressive time series model. *Biometrika* **68** 189–196. MR0614955

[22] HALL, P. and CARROLL, R. J. (1989). Variance function estimation in regression: The effect of estimating the mean. *J. Roy. Statist. Soc. Ser. B* **51** 3–14. MR0984989

[23] HALL, P. and TITTERINGTON, D. M. (1988). On confidence bands in nonparametric density estimation and regression. *J. Multivariate Anal.* **27** 228–254. MR0971184

[24] HÄRDLE, W. (1989). Asymptotic maximal deviation of $M$-smoothers. *J. Multivariate Anal.* **29** 163–179. MR1004333

[25] HÄRDLE, W. and MARRON, J. S. (1991). Bootstrap simultaneous error bars for nonparametric regression. *Ann. Statist.* **19** 778–796. MR1105844

[26] JOHNSTON, G. J. (1982). Probabilities of maximal deviations for nonparametric regression function estimates. *J. Multivariate Anal.* **12** 402–414. MR0666014

[27] KNAFL, G., SACKS, J. and YLVISAKER, D. (1985). Confidence bands for regression functions. *J. Am. Statist. Assoc.* **80** 683–691. MR0803261

[28] KOMLÓS, J., MAJOR, P. and TUSNÁDY, G. (1975). An approximation of partial sums of independent RV's and the sample DF. I. *Z. Wahrsch. Verw. Gebiete* **32** 111–131. MR0375412

[29] PRIESTLEY, M. B. (1988). *Nonlinear and Nonstationary Time Series Analysis*. Academic Press, London. MR0991969

[30] ROBINSON, P. M. (1997). Large-sample inference for nonparametric regression with dependent errors. *Ann. Statist.* **25** 2054–2083. MR1474083

[31] SUN, J. and LOADER, C. R. (1994). Simultaneous confidence bands for linear regression and smoothing. *Ann. Statist.* **22** 1328–1345. MR1311978

[32] TONG, H. (1990). *Nonlinear Time Series Analysis. A Dynamical System Approach*. Oxford Univ. Press. MR1079320

[33] VERHOEVEN, P., PILGRAM, B., MCALEER, M. and MEES A. (2002). Non-linear modelling and forecasting of S&P 500 volatility. *Math. Comput. Simulation* **59** 233–241. MR1926366

[34] WU, W. B. (2003). Empirical processes of long-memory sequences. *Bernoulli* **9** 809–831. MR2047687

[35] WU, W. B. (2005). Nonlinear system theory: Another look at dependence. *Proc. Natl. Acad. Sci. USA* **102** 14150–14154. MR2172215

[36] WU, W. B. (2007). Strong invariance principles for dependent random variables. *Ann. Probab.* **35** 2294–2320. MR2353389

[37] WU, W. B. and ZHAO, Z. (2007). Inference of trends in time series. *J. Roy. Statist. Soc. Ser. B* **69** 391–410. MR2323759

[38] XIA, Y. (1998). Bias-corrected confidence bands in nonparametric regression. *J. Roy. Statist. Soc. Ser. B* **60** 797–811. MR1649488

[39] ZHAO, Z. and WU, W. B. (2006). Kernel quantile regression for nonlinear stochastic models. Technical Report No. 572, Dept. Statistics, Univ. Chicago.

DEPARTMENT OF STATISTICS
PENNSYLVANIA STATE UNIVERSITY
UNIVERSITY PARK, PENNSYLVANIA 16802-2111
USA
E-MAIL: zuz13@stat.psu.edu

DEPARTMENT OF STATISTICS
UNIVERSITY OF CHICAGO
5734 S. UNIVERSITY AVENUE
CHICAGO, ILLINOIS 60637
USA
E-MAIL: wbwu@galton.uchicago.edu